\documentstyle[12pt]{amsart}

\newtheorem{thm}{Theorem}[section]
\newtheorem{prop}[thm]{Proposition}
\newtheorem{cor}[thm]{Corollary}
\newtheorem{defn}[thm]{Definition}
\newtheorem{lemma}[thm]{Lemma}

\input psfig

\begin{document}
\title[Teichm\"{u}ller spaces of Klein surfaces]{Some results on
Teichm\"uller spaces of Klein surfaces}
\author{Pablo Ar\'{e}s Gastesi}
\date{July 20, 1995}
\maketitle
\begin{abstract}
In this paper, we prove some isomorphisms theorems between Teichm\"{u}ller
spaces of
non-orientable compact surfaces.
We also develop a technique, based on similar results for Riemann
surfaces, to give explicit examples of Teichm\"{u}ller spaces of Klein
surfaces.
\end{abstract}

\section{Background and statement of main results}

The deformation theory of non-orientable surfaces
deals with the problem of studying parameter spaces for the different
dianalytic structures that a surface can have.
It is an extension of the classical theory
of Teichm\"{u}ller spaces of Riemann surfaces, and as such, it is quite rich.
In
this paper we study some basic properties of the Teichm\"{u}ller
spaces of non-orientable
surfaces, whose parallels in the orientable situation are well known.
More precisely, we prove an uniformization theorem, similar to the case
of Riemann surfaces, which shows that a non-orientable compact surface can
be represented as the quotient of a simply connected domain of the
Riemann sphere, by a discrete group of
M\"{o}bius and anti-M\"{o}bius transformation (mappings whose
conjugates are M\"{o}bius transformations).
This uniformization result allows us to give explicit
examples of Teichm\"{u}ller spaces of non-orientable surfaces,
as subsets of deformation spaces of orientable surfaces.
We also prove two isomorphism theorems: in the first place, we show that
the Teichm\"{u}ller spaces of surfaces of different topological type are not,
in general,
equivalent. We then show that, if the topological type is preserved, but the
signature changes, then the deformations spaces are isomorphic. These are
generalizations of the Patterson and Bers-Greenberg theorems for
Teichm\"{u}ller spaces
of Riemann surfaces, respectively.

A Riemann surface $(\Sigma,X)$ is a topological surface $\Sigma$ with a
complex structure $X$, that is, a covering of $\Sigma$ by charts with
holomorphic changes of coordinates. Since holomorphic functions have positive
Jacobian, it turns out that Riemann surfaces are orientable.
The natural generalisation
to the case of non-orientable surfaces is that of a {\bf dianalytic
structure},
where we require that the changes of coordinates are either holomorphic
or anti-holomorphic (the complex conjugate is holomorphic). A pair
$(\Sigma,X)$, where $\Sigma$ is a surface and $X$ is a dianalytic structure,
is called a {\bf Klein surface}. In particular, Riemann surfaces are Klein
surfaces. It is classical fact that any Klein surface can be
represented as $\tilde{X}/\Gamma$, where $\tilde X$ is either the Riemann
sphere,
the complex plane or the upper half plane, and $\Gamma$ is a group of
dianalytic bijections of $\tilde X$.  Except for a few (finite number of)
cases, Klein surfaces are covered by the upper half plane; these are called
{\bf hyperbolic surfaces}. A compact non-orientable
surface $\Sigma$ is
the connected sum of $g$ (real) projective planes; $g$ is called the
{\bf genus} of the surface. Observe that here we use the genus in the
topological sense; some authors (in particular, \cite{sep:book})
use the so-called arithmetic genus, which is equal to $g-1$.
A non-orientable surface is hyperbolic if and only if
$g\geq 3$. In the first result of this paper, we prove a uniformization
theorem, by groups which are more suitable for computations that
groups acting on the upper half plane.
\begin{thm}\label{thm-unif}
Let $\Sigma$ be a compact non-orientable surface of genus bigger
than $2$. Then
there exists a Kleinian group $G$, acting discontinuously on
a simply connected
set $\Delta$ of $\hat{\Bbb C}$, and an antiholomorphic function $r$, such
that:\\
1. $g(\Delta)=\Delta$ for all $g \in G$; $r(\Delta)=\Delta$; \\
2. $\Delta/G$ is isomorphic to the complex double $\Sigma^c$ of $\Sigma$;\\
3. $r$ is of the form $r:z \rightarrow \frac{a\overline{z}+b}
{c\overline{z}+d}$, with $ad-bc\neq 0$;\\
4. $\Delta/\Gamma\cong \Sigma$, where $\Gamma$ is the group generated by
$G$ and $r$;\\
5. $\Gamma$ is unique up to conjugation by M\"{o}bius transformations.\end{thm}
Here by a Kleinian group we mean a group of M\"{o}bius
transformations that acts discontnuously on a non-empty open set of the
Riemann sphere.
The {\bf complex double} of $\Sigma$ is a Riemann surface $\Sigma^c$,
together with a unramified double cover $\pi:\Sigma^c \rightarrow \Sigma$.
If $\Sigma$ is hyperbolic, then $\Sigma^c$ is also hyperbolic (see \S 2 below).

Let $M(\Sigma)$ denote the set of dianalytic structures, on the non-orientable
surface $\Sigma$, that are compatible with the differential structure
induced by $X$. The quotient of $M(\Sigma)$ by the group of diffeomorphisms
homotopic to the identity (acting by pullback, see \S $3$),
is the {\bf Teichm\"{u}ller space} $T(\Sigma)$ of $\Sigma$. It has a natural
real analytic
structure given by projecting the natural structure of $M(\Sigma)$.
It is not hard to prove that $T(\Sigma)$ embedds in the Teichm\"{u}ller space
of $\Sigma^c$ (see \S $3$). Combining this embedding with theorem
\ref{thm-unif} and the results of I. Kra in \cite{kra:horoc},
we can give presentations for the deformation
spaces of some non-orientable surfaces. As an example, we compute the
Teichm\"{u}ller space of a surface of genus $3$.
\begin{thm}\label{thm-ex1}
The space $T(\Sigma)$ of a non-orientable surface of genus $3$
can be identified with the set of points $(\tau_1,\tau_2,\tau_3)$ of
$T(\Sigma^c)$, such that
$$\left\{\begin{array}{lcl}
{\rm Re}(\tau_2)		  & = & 0\\
{\rm Re}(\tau_1)		  & = & {\rm Im}(\tau_3)\\
{\rm Re}(\tau_1)+{\rm Re}(\tau_3) & = & 0
\end{array}\right .$$
\end{thm}
We introduce the concept of puncture on a non-orientable surface as a
generalisation of the corresponding idea on Riemann surfaces: a puncture
is a domain on $\Sigma$, homeomorphic to the unit disc minus the origin,
that cannot be completed to be homeomorphic to the unit disc, and such
that any change of coordinates in the domain is holomorphic.
The above theorems extends easily to the case of surfaces
with punctures. For example, we can identify the deformation space of a
surface of genus $1$ with two punctures.
\begin{thm}\label{thm-ex2}
The space $T(\Sigma)$, where $\Sigma$ is the (real) projective
plane with two puntures, can be identified with the set of points of
the upper half plane with imaginary part bigger then $1$.\end{thm}
One can define a {\bf Klein hyperbolic orbifold} as a non-orientable surface
$\Sigma$, with finitely many (maybe zero) punctures,
such that the covering from the upper half plane to $\Sigma$ is ramified
over a finiter number of points.
The surface $\Sigma^c$ carries an anticonformal involution $\sigma$,
such that $\Sigma^c/<\sigma> \cong \Sigma$. We have that $T(\Sigma)$ can
be identified with the set of fixed points of the anticonformal involution
$\sigma^*$, induced by $\sigma$ in $T(\Sigma)$. We say that the Teichm\"{u}ller
spaces of
two non-orientable surface $\Sigma_1$ and $\Sigma_2$ are {\bf real isomorphic},
if there exists a biholomorphic mapping $f: T(\Sigma_1^c) \rightarrow
T(\Sigma_2^c)$, such that $f\circ \sigma_1^* = \sigma_2^* \circ f$.
The following result is a generalisation of the Bers-Greenberg
isomorphism for Riemann surfaces.
\begin{thm} [Bers-Greenberg theorem for non-orientable surfaces]
If $\Sigma_i$,
$i=1,2$, are two non-orientable hyperbolic orbifolds, with the
same genus and number of ramification points, then the spaces
$T(\Sigma_1)$
and $T(\Sigma_2)$ are real isomorphic.\label{thm-bers}\end{thm}

This paper is organized as follows: in \S 2 we prove the uniformization
theorem; \S 3 contains the proof of theorem \ref{thm-bers} and other
results about isomorphisms of deformation
spaces; finally, \S 4 has the examples.

{\bf Akcnowledgements}: I would like to thank R. R. Simha, for useful
discussions. I would like also to thank the Tata Institute of
Fundamental Research, for providing me with excellent research
facilities.

\section{Uniformization}

Classically, hyperbolic Klein surfaces are uniformized as the
quotient of the upper half plane by a discrete group of dianalytic
self-homeomorphisms (M\"{o}bius and anti-M\"{o}bius transformations with real
coefficients),
known as NEC groups.
In this section, we will prove a uniformisation theorem by a different
type of groups, which are more suitable for computations. We will use these
groups, in \S $4$, to produce some {\it explicit} examples of deformation
spaces of non-orientable surfaces.

We start by recalling some facts of uniformization of Riemann surfaces. A
{\bf partition} $\cal C$ on a Riemann surface $\Sigma$, of genus $g \geq
2$, is a collection of simple closed disjoint curves, such that no
curve of $\cal C$ is homotopically trivial, and no two curves of $\cal C$
are freely homotopically equivalent. A partition consists of at most
$3g-3$ curves; if this bound is attained, we say that the partitions is
{\bf maximal}. See \cite{str:qd} for the proof of the existence of
partition on surfaces.
\begin{thm}[Maskit Uniformization Theorem, \cite{mas:bound}, \cite{mas:class}]
Given a \\ Riemann surface $\Sigma$
and a maximal partition ${\cal C} = \{a_1,\ldots,a_{3g-3+n}\}$, there
exists a Kleinian group $G$, known as a terminal regular b-group, such
that:\\
1. there is a unique maximal simply connected set $\Delta$ of the Riemann
sphere, where $G$ acts discontinuously , and $g(\Delta)=\Delta$
for all $g \in G$;\\
2. $\Delta/\Gamma \cong \Sigma$;\\
3. to each curve of $\cal C$ corresponds a maximal conjugacy class of
cyclic subgroups of $G$ generated by a parabolic transformation;\\
4. besides $\Sigma$, the group $G$ uniformizes the $2g-2$ thrice
punctured spheres obtained from squeezing each curve of $\cal C$ to a
puncture;\\
5. $G$ is unique up to conjugation by M\"{o}bius transformations.
\end{thm}
A {\bf symmetry} $\sigma$ on a Riemann surface is an anticonformal involution.
If $F(\sigma)$ denotes the set of fixed points of $\sigma$, then we have
that $\Sigma - F(\sigma)$ consists of at most two components. It is a
well known fact that $\Sigma/<\sigma>$ is orientable if and only if
$\Sigma - F(\sigma)$ is not connected (\cite{buja:symm}, \cite{sep:book}).
The classical result about the structure of $F(\sigma)$ is the following.
\begin{thm}[Harnack] If $\sigma$ is a symmetry on a compact surface $\Sigma$
of genus $g$, then $F(\sigma)$ is either empty or
consists of $s$ simple disjoint curves
$\delta_j$, with $s\leq g+1$.\end{thm}
\noindent This theorem can be improved as follows.
\begin{thm}[Kra-Maskit] In addition to the curves (if any) $\delta_1,
\ldots, \delta_s$, there exists closed curves $\delta_{s+1}, \ldots, \delta_t$,
such that:\\
1. $\{\delta\}_{j=1}^t$ is a collection of disjoint curves;\\
2. $\sigma(\delta_j)=\delta_j$, for all $j$;\\
3. $\Sigma - \cup_{j=1}^t\delta_j$ consists of two components, $\Sigma_1$
and $\Sigma_2$;\\
4. $\sigma$ interchanges $\Sigma_1$ and $\Sigma_2$.\end{thm}
The existence of maximal partitions invariant under symmetries
is a well known fact; but our proof is different from those in the
literature (see, for example \cite[pgs. 117-120]{sep:book}),
but we include it here for the sake of completeness.
\begin{lemma}\label{lemma-part}
Let $\Sigma$ be a Riemann surface of genus $g \geq 2$, and let
$\sigma$ be a symmetry on $\Sigma$. Then there exists a maximal partition
$\cal C$ invariant under $\sigma$, that is, $\sigma({\cal C}) = \cal C$.
\end{lemma}
\begin{pf}Let ${\cal C}_1$ denote the set of curves given by the
Harnack-Kra-Maskit theorems. We claim that ${\cal C}_1$ is a partition
on $\Sigma$. In fact, we have that if a curve of ${\cal C}_1$ is
homotopically trivial, then $\Sigma_1$ and $\Sigma_2$ are discs, and
therefore, $\Sigma$ will be homeomorphic to the Riemann sphere. Similarly,
if two curves of ${\cal C}_1$ are freely homotopic, we get that $\Sigma$
is a torus.

If ${\cal C}_1$ is maximal, we are done. If not, let $a$ be a curve such
that ${\cal C}_2 = {\cal C}_1 \cup \{a\}$ is a partition. We claim
that ${\cal C}_2 \cup \{\sigma(a)\}$ is a partition. This can be seen in
three easy steps:\\
1. $\sigma(a)$ is not homotopically trivial, since $\sigma$ is a homeomorphism,
and $a$ is not trivial (being a curve in a partition);\\
2. $\sigma(a)$ is not (freely) homotopically equivalent to any curve
of ${\cal C}_1$. If there is a curve $\delta$ in ${\cal C}_1$, freely
homotopic to $\sigma(a)$, then, applying $\sigma$, we would get that
$a$ is freely homotopic to $\delta$, contradicting the fact that ${\cal C}_2$
is a partition;\\
3. $\sigma(a)$ is not freely homotopic to $a$. If these two curves are
freely homotopic, then
we have that $a$ and $\sigma(a)$ bound a cylinder in $\Sigma$. Since
these curves lie in different components of $\Sigma - {\cal C}_1$, we
get that there is a curve, $\delta$ in ${\cal C}_1$,
in that cylinder. But this implies that $a$ is homotopically equivalent
to $\delta$, which is again not possible.
\end{pf}
Any non-orientable surface $\Sigma$ has a double unramified cover by a
Riemann surface $\Sigma^c$, called the {\bf complex double}
(\cite[pgs 37-40]{all:klein}). If $\Sigma$ has genus $g$, then
$\Sigma^c$ has genus $g-1$. $\Sigma$ has a symmetry
$\sigma$, such that $\Sigma^c/<\sigma> \cong \Sigma$.
We have now all the necessary tools to prove theorem \ref{thm-unif}.
\begin{pf*}{Proof of theorem \ref{thm-unif}}Let $\Sigma^c$ be the complex
cover of
$\Sigma$, and let $\sigma$ be the symmetry on $\Sigma^c$ such that
$\Sigma^c/<\sigma> \cong \Sigma$. By our hypothesis, $\Sigma^c$ has
genus greater than $1$, so applying the lemma \ref{lemma-part} we
obtain a $\sigma$-invariant maximal partition
$\cal C$ on $\Sigma^c$. Using the Maskit Uniformization
Theorem, we get a Kleinian group $G$, uniformizing $\Sigma^c$ in the invariant
simply connected component $\Delta$. We only need to show
that the symmetry $\sigma$ lifts to an anti-M\"{o}bius transformation,
in the covering determined by $G$
(i.e., it is of the form given in the statement of the theorem).
For simplicity, assume first that $\sigma$ is orientation preserving.
Then $\sigma$ induces a set of conformal mappings, $\sigma_j:S_j
\rightarrow S_k$, among the parts $S_1, \ldots, S_{2g-2}$ of $\Sigma - \cal C$.
The infinite Nielsen extension, $\tilde{S_j}$, of $S_j$ is a thrice punctured
sphere, obtained from $S_j$ by completing the holes to punctured discs. It
is a classical fact that $\sigma_j$ extends to a quasiconformal mapping,
denoted by $\tilde{\sigma_j}$, from $\tilde{S_j}$ to $\tilde{S_k}$, with
maximal dilatation $1\leq K(\tilde{\sigma_j}) \leq K(\sigma_j)$
(\cite{bers:nielsen}, \cite{abik:book}). Since $\sigma_j$ is conformal, we
have that its
dilatation is equal to $1$, and therefore $K(\tilde{\sigma_j})=1$, that
is, $\tilde{\sigma_j}$ is also conformal.  Let $\Delta_j$
be a component of $\pi^{-1}(S_j)$, where $\pi:\Delta \rightarrow \Sigma$
is the natural quotient mapping from $\Delta$ onto $\Sigma$, and let
$G_j=stab(G,U_j):=\{g\in G;~g(U_j)=U_j \}$. We have that the $G_j'$s
are triangle
groups with two invariant components; let $U_j'$ be the component that does
not contain $\Delta$. Then, the mapping $\sigma_j$ induces a conformal
mapping between $U_j'$ and $U_k'$, for a proper choice of $U_k$. This can
be done with all the components of $\pi^{-1}(S_j)$, and all the $j=1,
\ldots, 2g-2$, obtaining in this way a conformal self-mapping
$\tilde{\sigma}$, of $\Delta
\cup g(U-j)$. But this set is the region of discontinuity of $G$
(that is, the set of points of the Riemann sphere were $G$ acts
discontinuously). Since $G$ is finitely generated, we have that the complement
of the region of discontinuity has measure zero. Therefore, the classical
theory of quasiconformal mappings gives us a conformal automorphism of
the Riemann sphere that extends
$\tilde{\sigma}$. Such mapping should be a M\"{o}bius transformation.

To complete the proof of the theorem it suffices to observe the following two
facts, which are easy to prove:\\
1. the theory of quasiconformal mappings has a natural extension
to cover the orientation reversing mappings \cite{sep:book}; and\\
2. Bers' results on Nielsen extensions can be applied to orientation
reversing mappings.\end{pf*}
We define a {\bf ramification point} $x$ on a Klein surface as a point
such that the universal covering looks like $z \mapsto z^n$, in a
neighborhood of $x$, (which corresponds to the points $z=0$) for some finite
positive integer $n$.
The number $n$ is called the {\bf ramification value} of $x$. Ramification
points correspond to fixed points of orientation preserving
transformations, of finite order.
\begin{defn}A {\bf puncture} is a domain
$D$ in $\Sigma$ satisfying the following conditions:\\
1. $D$ is homeomorphic to ${\Bbb D}^*=\{z\in {\Bbb C};~0<|z|<1\}$;\\
2. for any sequence of points in ${\Bbb D}^*$ converging to the origin,
the corresponding sequence in the surface diverges;\\
3. if there are two patches on $\Sigma$, whose images contain some sets
of the form $\{z\in {\Bbb C};~0<|z|<r\leq 1\}$ (that is, neighborhoods
of the \lq\lq missing point\rq\rq),
then the change of coordinates is holomorpic.\end{defn}

Given a Klein surface with ramification points and/or
punctures, called a {\bf Klein orbifold},
we define its {\bf signature} as a collection of numbers (and
a symbol) of the form $(g, \pm, n; \nu_1, \ldots, \nu_n)$, where
$g$ is the genus of the surface, $n$ is the number of special points, and
$\nu_1, \ldots, \nu_n$ are the ramification values, with puntures having
ramification value equal to $\infty$. If the orbifold is orientable, then
we take the symbol $+$, while $-$ is used for non-orientable surfaces.
If all the ramification values are equal to $\infty$, then we will write
the signature as $(g, \pm, n)$.
It is not difficult to see that if $\Sigma$ has signature $(g, -, n;
\nu_1, \ldots, \nu_n)$, then the signature of $\Sigma^c$ must be $(g-1,
+, 2n; \nu_1, \nu_1, \ldots, \nu_n, \nu_n)$.
A Klein orbifold $\Sigma,$ is hyperbolic if and only if
$kg-2+n-\sum_{j=1}^n\frac{1}{n}$ is positive, where $k=1$ if $\Sigma$ is
not orientable, and $k=2$ in the orientable case.
Since the Maskit Uniformization Theorem and the theory of quasiconformal
mappings extend to the case of surfaces with ramification points,
we have that theorem \ref{thm-unif} can be applied also to hyperbolic Klein
orbifolds.

\section{Isomorphisms between Teichm\"{u}ller spaces}

A natural problem in deformation theory is to study
which properties of a surface are determined by its Teichm\"{u}ller space, and
vice
versa.
More precisely, in this section we will see that, if the Teichm\"{u}ller spaces
of
two surfaces are equivalent, then the surfaces are homeomorphic.
Reciprocally, if two Klein orbifolds have the same genus and
number of ramification points, we will prove that their deformation
spaces are isomorphic.

We start by recalling the definition of the modular group, and
some basic facts about hyperelliptic surfaces.
The {\bf modular group} $Mod(\Sigma)$ of a non-orientable surface is
the quotient of the group of diffeomorphisms, by those homotopic to
the identity (in the case of Riemann surfaces, one takes only the
orientation preserving diffeomorphisms).
We have that $Mod(\Sigma)$ acts on
$T(\Sigma)$ by pullback: given a mapping $f$, and a real analytic structure
structure $X$, we define $f^*(X)$ as the unique structure on $\Sigma$
that makes $f:(\Sigma,f^*(X)) \rightarrow (\Sigma,X)$ dianalytic.
The mapping $f^*:[X] \rightarrow [f^*(X)]$
becomes dianalytic in the natural structure of $T(\Sigma)$.

We say that a non-orientable surface $\Sigma$ is {\bf hyperelliptic} if it
is a double cover of the (real) projective plane (respectively, the
Riemann sphere, in case of orientable surfaces).
Hyperelliptic surfaces carry the so-called {\bf hyperelliptic
involution}, which is a dianalytic involution (holomorphic, in the
case of orientable surfaces) $\alpha$, such that $\Sigma/<\alpha>$ is the
projective plane. It is not hard to see that if $\Sigma$ is hyperelliptic,
so is its complex cover $\Sigma^c$, and that $\alpha$ lifts to the
hyperelliptic involution $j$ of $\Sigma^c$.
Since the hyperelliptic involution on a
Riemann surface is unique, we obtain the reciprocal result: if
$\Sigma^c$ is the complex double of a surface $\Sigma$, and $\Sigma^c$
is hyperelliptic, then $\Sigma$ is also hyperelliptic; moreover, the
involution $j$ can be pushed down to the hyperelliptic involution
$\alpha$ in $\Sigma$.

As in the case of Riemann surfaces, we have that the modular group acts
effectively on Teichm\"{u}ller space, except for a finite number of cases.
\begin{prop} $Mod(\Sigma)$, for a hyperbolic non-orientable surface
$\Sigma$, compact with finitely (maybe zero) punctures, acts effectively on
$T(\sigma)$, except for the following cases:
the projective plane with two puntures, the Klein bottle with one
puncture, or the connected sum of three projective planes.\end{prop}
\begin{pf}As usual, let $\Sigma^c$ denote the
complex double of $\Sigma$, and let $\sigma$ be the involution
associated to such covering. It is a well known fact (\cite[pg.
149]{sep:book} or below) that $T(\Sigma)$ can be identified with the set of
fixed points $T(\Sigma^c)_{\sigma^*}$ of $\sigma^*$ in $T(\Sigma^c)$.
If $f$ is a diffeomorphism of $\Sigma$, there is a unique
holomorphic lift, $F:\Sigma^c \rightarrow \Sigma^c$ (see
\cite[pg. 20]{sep:spaces} and \cite[pg. 39]{all:klein}).
By uniqueness, we have that $F$
satisfies $F^*\circ \sigma^* = \sigma^*\circ F^*$.
This proves that $Mod(\Sigma)$ embeddes into
the set $A=\{h^* \in Mod(\Sigma^c);~ h^*\circ \sigma^* = \sigma^*\circ
h^*\}$.
The modular group of a hyperbolic riemann surface acts properly on the
corresponding Teichm\"{u}ller space, unless the signature of the surface is
$(0, +, 4)$, $(1, +, 2)$, $(2, +, 0)$ or $(1, +, 1)$. Since complex
covers have an even number of punctures, we get that only the first
three signatures can give non-orientable surfaces. We therefore obtain
that,
the cases where $Mod(\sigma)$ may fail to act effectively correspond to
the signatures $(1, -, 2)$, $(2, -, 1)$ and $(3, -, 0)$.
This proves the first part of the proposition.

If, in the other hand, $\Sigma$ has signature in the above list, we have
that the only elements of $Mod(\Sigma)$ that do not act properly are the
classes of the identity and the hyperelliptic involution.
By the remarks before the theorem, we also have that these classes are
the only elements that act trivially on the Teichm\"{u}ller space of hyperbolic
Klein
surfaces. \end{pf}
It would be interesting to know whether the image of $Mod(\Sigma)$ is
equal to the whole set $A$ described in the above proof.

The proof of the following result is straightforward from the parallel
result for Riemann surfaces. Nevertheless, the proposition is
interesting, because it shows the great similarity between the theory
of deformation of orientable surfaces, and that of non-orientable ones.
\begin{thm}
[Patterson theorem for non-orientable surfaces]
Let $\Sigma_i$, $i=1,2$, be two hyperbolic compact, with finitely many
(possibly zero) punctures, non-orientable surfaces, and suppose that
either $\Sigma_1$ or $\Sigma_2$ has genus not equal to $3$.
If $T(\Sigma_1)$ is real isomorphic to
$T(\Sigma_2)$, then $\Sigma_1$ is homeomorphic to $\Sigma_2$.\end{thm}
\begin{pf}It suffices to observe that if $T(\Sigma_1)$ and $T(\Sigma_2)$
are real isomorphic, then the spaces $T(\Sigma_1^c)$ and $T(\Sigma_2^c)$
are biholomorphic, and therefore, $\Sigma_1^c$ and $\Sigma_2^c$ are
homeomorphic.\label{thm-pat}\end{pf}
In order to prove theorem \ref{thm-bers}, we need to review a basic
concept of Teichm\"{u}ller theory: quadratic differentials and Beltrami
coefficients.
Let $Q(\Sigma)$ denote the space of {\bf bounded quadratic differentials} on
a Klein surface. These are simply quadratic differentials, regular on the
surface, with at most simple poles at the punctures, and with zeros
of certain order (determined by the ramification value) at the
ramification points (see, for example \cite{kra:autom}). It is easy
to see that the set $Q(\Sigma)$, for a non-orientable surface, can be
identified with the subspace of elements of $Q(\Sigma^c)$, that are preserved
by $\sigma$, that is, $\overline{\phi} =
(\phi\circ \sigma) (\overline{\partial}\sigma)^2$.
Here by $\overline{\partial}\sigma$ we mean
$\partial\sigma/\partial\overline{z}=\frac{1}{2}(\frac{\partial}{\partial
x}+i \frac{\partial}{\partial y}) \sigma$.
The dimension of $Q(\Sigma)$ over ${\Bbb R}$
is equal to the dimension of $Q(\Sigma^c)$ over $\Bbb C$.
The embedding from $Q(\Sigma)$ into $Q(\Sigma^c)$ is an isometry in the
norm, $$||\varphi|| = sup~|\varphi(x)|\lambda^{-2}(x),\hspace{5mm}f \in
Q(\Sigma),$$ where $\lambda$ is
the metric obtained by pushing down the Poincar\'{e} metric of the upper
half plane onto the corresponding surface (and the supremum is taking
over the whole surface).
The cotangent bundle of $T(\Sigma)$ can be naturally identified with the
space $Q(\Sigma)$.

Let $\Sigma$ be a Klein surface, uniformized by the NEC group $\Gamma$.
Let $M({\Bbb H},\Gamma)$ denote the space of {\bf Beltrami differentials}
for $\Gamma$. This set consists of (classes of)
measurable functions $\mu$, with support in the upper half plane, and
$L^\infty$-norm less than one,
satisfying $(\mu\circ\gamma)\overline{\gamma'}/\gamma
= \mu$, if $\gamma \in \Gamma$ is orientation preserving, or
$(\mu\circ\gamma)\overline{\overline{\partial}\gamma'}/
\overline{\partial}\gamma = \overline{\mu}$,
if $\gamma \in \Gamma$ reverses the orientation. For each $\mu \in
M({\Bbb H},\Gamma)$, there is a unique quasiconformal homeomorphism $w_\mu$,
of the upper half plane, with dilatation $\mu$, that fixes $0$, $1$ and
$\infty$. Two Beltrami coefficients, $\mu$ and $\nu$, are {\bf equivalent}
if $w_\mu = w_\nu$ on the real line. The space of Beltrami differentials,
quotiented by the above equivalence relation is the Teichm\"{u}ller space
$T(\Gamma)$
of the group $\Gamma$. It can be proven that $T(\Gamma)$ is naturally
isomorphic to
$T(\Sigma)$, where $\Sigma \cong {\Bbb H}/\Gamma$.
It is easy to see that, if $G$ is the subgroup of $\Gamma$
consisting of the orientation preserving mappings, then ${\Bbb H}/G \cong
\Sigma^c$. We can identify the Beltrami differentials for $\Gamma$
with those Beltrami differentials for $G$, that are invariant under
$\sigma$, that is,
$(\mu\circ\sigma)\overline{\overline{\partial}\sigma}/
\overline{\partial}\sigma = \mu$.
This allows us to identify the deformation space $T(\Sigma)$ with the
set of fixed points of $\sigma^*$,
$T(\Sigma^c)_{\sigma^*}$, in the deformation space of
the complex double $T(\Sigma^c)$.
In this way we obtain a Teichm\"{u}ller's lemma for non-orientable surfaces:
on each equivalence class of Beltrami differentials there is a unique mapping
with minimal dilatation, which is of the form $\mu=k\overline{\varphi}/
|\varphi|$, with $\varphi\in Q(\Sigma)$, $k$ a real number.
With this background, we can provide two proofs of theorem \ref{thm-bers} of
\S $1$.
\begin{pf*}{First proof of theorem \ref{thm-bers}}
Let $\Sigma$
be a hyperbolic surface of signature $(g,-,n;$
$\nu_1,\ldots, \nu_n)$, where
we assume that at least one of the ramification values is finite. Let
$\Sigma_0$ be
the surface of signature $(g, -, n; \infty, \ldots, \infty)$, obtained
by removing from $\Sigma$ all the points with finite ramification value.
Let $\Gamma$ and $\Gamma_0$ be NEC groups uniformizing $\Sigma$ and
$\Sigma_0$ respectively.
Define ${\Bbb H}_{\Gamma}$ as ${\Bbb H} - \{$ fixed points of elliptic
elements of $\Gamma \}$. Then, by our hypothesis we have that ${\Bbb H}
\neq {\Bbb H}_\Gamma$. Since $\Sigma_0 \cong {\Bbb H}_\Gamma/\Gamma$,
we have a covering map $h:{\Bbb H} \rightarrow {\Bbb H}_\Gamma$, that makes
diagram $1$ commutative.
The function $h$ induces a group homomorphism $\chi:\Gamma_0
rightarrow
\Gamma$, defined by the rule $h\circ \chi(\gamma) = \gamma\circ h$.
The mapping $h$ induces a mapping, $h^*:M({\Bbb H},\Gamma_0) \rightarrow
M({\Bbb H},\Gamma)$, between Beltrami coefficients,
given by the expression $(h^*\mu)\circ h = \mu h'/\overline{h'}$ (see
below for a proof of the fact that $h$ is a holomorphic function).
It is not hard to see, using the same arguments that in the orientable case,
that $h^*$ induces an real analytic bijection between the
spaces $T(\Sigma_0)$ and $T(\Sigma)$. See, for example, \cite{ek:hol}
(or \cite{ares:bers}, for a more details).
By analytic continuation, we can extend $h^*$ to a biholomorphic
function between $T(\Sigma_0^c$ and $T(\Sigma^c$, that obviously
commutes with the involutions $\sigma_0$ and $\sigma$, giving the
desired isomorphism.
\end{pf*}
\begin{pf*}{Second proof of theorem \ref{thm-bers}}In this case, we will use
the Bers-Greenberg theorem for Riemann surfaces. Consider the same setting
as in the first proof. Let $G$ and $G_0$ be the orientation preserving
subgroups of $\Gamma$ and $\Gamma_0$ respectively. Then we have ${\Bbb
H}/G \cong \Sigma^c$ and ${\Bbb H}/G_0 \cong {\Bbb H}_G/G \cong
\Sigma_0^c$. We get that diagram $2$ is commutative.
The functions $\pi$, $\pi_0$ and $\rho$ are the natural projections. The
mapping $h$ is defined as in the first proof. The function $f$ is the
unique holomorphic mapping that makes the lower triangle commutative. We
have that $\pi$ is holomorphic (a covering of a Riemann surface by an open
set of the complex plane), so the function $h$ is holomorphic. This
implies that the group homomorphism $\chi$, of the first proof, takes the
subgroup $G_0$ onto $G$. Since the surfaces $\Sigma$ and $\Sigma^c$ have
the same genus and number of ramification points/punctures, we have that
the spaces $T(\Sigma^c_0)$ and $T(\Sigma^c)$ are isomorphic, via the
function $h^*$ induced by $h$, as in the first proof.
To prove theorem \ref{thm-bers} it suffices to show that
$h$
commutes with the antiholomorphic involutions $\sigma$ and $\sigma_0$
(that \lq\lq produce\rq\rq\ the surfaces $\Sigma$ and $\Sigma_0$,
respectively). In other words, we have to show $h\circ \sigma_0 =
\sigma\circ h$. But we have that $\pi\circ h\circ\sigma_0 =
\pi\circ \sigma\circ h$, so $h\circ\sigma_0$ is equal to either
$\sigma\circ h$ or $\sigma\circ h\circ \sigma$. Since this last function
is holomorphic, we must have $h\circ \sigma_0 = \sigma\circ h$,
as claimed. Identifying
$T(\Sigma_0)$ and $T(\Sigma)$ with the set of fixed points of $\sigma_0^*$
and $\sigma^*$ in $T(\Sigma_0^c)$ and $T(\Sigma^c)$, respectively, we get
the Bers-Greenberg theorem for non-orientable surfaces.
\end{pf*}

\section{Examples}

In this section, we will shown with two examples, how the
techniques of Kra of \cite{kra:horoc} can be applied to the case of
non orientable surfaces. We will work with deformation spaces of
Kleinian groups, which are equivalent (if the groups are chosen
properly, for example,  groups given by theorem $1$) to deformation
spaces of Riemann or Klein surfaces (see \cite{kra:variational} or
\cite{sep:book} for more details).

In out first example, we consider a Klein
surface, $\Sigma$, of signature $(0,-,2)$. Its complex double,
$\Sigma ^c$,
has signature $(0,+,4)$. A Kleinian group, $G_\alpha$,
uniformizing $\Sigma^c$ is generated by the transformations
$$A=\left[\begin{array}{cc}-1&-2\\0&-1\end{array}\right],~
B=\left[\begin{array}{cc}-1&0\\2&-1\end{array}\right],~
B_\alpha= \left[ \begin{array}{cc} -1+2\alpha&-2\alpha^2\\2&
-1-2\alpha\end{array}\right],$$
where ${\rm Im}(\alpha)>1$ (see above reference).
The coordinate of $G_\alpha$ in the Teichm\"{u}ller space $T(0,+,4)$
(notation should be obvious)
is given by the expression
$$\alpha=cr(f(A),f(B),f(AB),f(B_\alpha)).$$
Here $cr$ denotes the cross ratio of
four points in the Riemann sphere, chosen so that $cr(\infty,0,1,z)=z$,
and $f(T)$ denotes the unique fixed point of the parabolic transformation
$T$.

A maximal partition in $\Sigma^c$ consists of a simple
closed curve, say $a_1$. We can assume that
the punctures $P_1$ and $P_2$ lie on the same component of
$\Sigma^c-a_1$. The M\"{o}bius transformation $A$ corresponds to the partition
curve.
Let $\gamma_j$, $j=1,\dots,4$, be a small simple loop around the
puncture $P_j$, oriented such that the puncture lies to the left of
$\gamma_j$.
The parabolic elements $T_1=B$, $T_2=(AB)^{-1}$, $T_3=B_\alpha^{-1}$
and $T_4=B_\alpha A$, correspond to these four loops.
Without loos of generality, we can assume that $T_j$ corresponds to
$\gamma_j$.
Consider on $\Sigma ^c$ the involution $\sigma=r\circ R$, where
$R$ is a rotation by $180$ degrees on the axis
of figure $3$, and $r$ is an anticonformal reflection on $a_1$.\vspace{2mm}
The anticonformal mapping
$\sigma$ has not fixed points, and the quotient $\Sigma ^c / < \sigma >$
has signature $(0, -, 2)$.

In order to identify $T(\Sigma)$ in $T(\Sigma^c)$, we have to study the
action of the mappings $R$ and $r$ on the group $G_\alpha$.
The transformation $R$ intechanges the punctures that lie on the same
component of $\Sigma^c-\{a_1\}$, that is, $R$ sends $\gamma_1$ to
$\gamma_2$ and $\gamma_3$ to $\gamma_4$ (up to free homotopy).
The function $r$, not only intechanges the two components of
$\Sigma^c-\{a_1\}$, but also changes the orientation of the loops,
sending $\gamma_1$ to $\gamma_4^{-1}$ and $\gamma_2$ to
$\gamma_3^{-1}$.
We have that $R$ lifts to $A_1^{1/2}$, while $r$ lifts to
$\tilde{r}(z)= z+\mu$, $\mu\in\Bbb C$, in the covering determined by
the group $G_\alpha$. Observe that, although a M\"{o}bius transformation may
have may square roots, parabolic elements have only one, and therefore,
$A_1^{1/2}$ is well defined.

{}From these obervations we can compute the action of $\sigma^*$ on
$T(\Gamma_\alpha)$ as follows. First of all, observe that the group
$G_\alpha$ is generated by $A$, $B$ and $B_\alpha$, with the
property that $AB$ and $A^{-1}B_\alpha^{-1}$ are parabolic
elements. We will use the notation $G(A,B,B_\alpha)$
to emphasis this fact. Observe also that $AB_\alpha^{-1}=B_{\alpha-1}$.
The mapping $R$ sends $G(A,B,B_\alpha)$ to
$G'=G(A,A^{1/2}BA^{-1/2},A^{1/2}B_\alpha A^{-1/2})=
G(A,B^{-1}A^{-1},B_{\alpha+1}^{-1})$.
Since the transformation $r$ is orientation reversing, we have that its
action is given by conjugating the group $G'$ into
$G(A,\tilde{r}B_{\alpha+1}\tilde{r}^{-1},$\\ $\tilde{r}AB\tilde{r}^{-1})$.
The mapping $\sigma^*$ has therefore the form
$$\sigma^*(\alpha)=cr(\infty,\overline{\alpha}+1+\mu,
\overline{\alpha}+2+\mu,1+\mu)=-\overline{\alpha}.$$
Therefore, the Teichm\"{u}ller
space $T(\Sigma)$ can be identified with the set of points $\alpha\in
T(\Sigma^c)$ such that ${\rm Re}(\alpha)=0$.
By the work of Kra, we have that $T(\Sigma)$ is precisely the set $\{ z
\in {\Bbb C};~{\rm Re}(z)=0,~{\rm Im}(z)>1 \}$.

Consider now the case of a surface $\Sigma$ of signature $(3,-,0)$.
The complex double of $\Sigma$
is a Riemann surface $\Sigma^c$, of genus $2$ without punctures. In
\cite{kra:horoc} we can find a group $G_\tau$ uniformizing
$\Sigma^c$,  generated by the M\"{o}bius transformations:
$$A_1=\left[\begin{array}{cc} -1&-2\\0&-1\end{array}\right],~
A_2=\left[\begin{array}{cc} 1&-2\\2&-3\end{array}\right],~
C_1=i\left[\begin{array}{cc} \tau_1&1\\1&0\end{array}\right],$$
$$A_3=\left[\begin{array}{cc} -1-2\tau_2 (1-\tau_2)&
-2(1-\tau_2)^2\\2\tau_2^2&-1+2\tau_2(1-\tau_2)\end{array}\right],$$
$$C_3=\left[\begin{array}{cc} \tau_3\tau_2^2 +2(1-\tau_3)\tau_2 +\tau_3-2
&-\tau_3\tau_2^2 + (3\tau_3-2)\tau_2-2\tau_3+3\\
\tau_3\tau_2^2 + (2-\tau_3)\tau_2-1&
-\tau_3\tau_2^2 + 2(1-\tau_3)\tau_2+2\end{array}\right].$$
The mapping $A_j$ correspond to the curves $a_j$ of the
partition of $\Sigma^c$ of figure $4$.
The $C_j$ are loxodromic
elements with the property that
$B_1^{-1}:=C_1^{-1} A_1 C_1$ and
$B_3:=C_3^{-1} A_3 C_3$ are parabolic elements. The coordinates on $T
(\Sigma^c)$ are given by the cross ratios
$$\left\{\begin{array}{lcl}
\tau_1=cr( f(A_1), f(B_1), f(A_2), C_1(f(A_1))),\\
\tau_2=cr( f(A_2), f(A_1), f(B_1), f(A_3)),\\
\tau_3=cr( f(A_3^{-1}), f(A_3A_2), f(A_2), C_3(f(A_3^{-1}))).
\end{array} \right .$$
Let $\sigma =rR$ be defined in a similar way as in the previous
example: $R$ is a rotation by $180$ degrees on the line of
figure $4$, and $r$ is an antiholomorphic reflection on the curve $a_2$.
The computation of the case $(0,-,2)$ applies to the coordinate
$\tau_2$, since the part corresponding to it (that is, $\Sigma
-\{a_1,a_3\}$) is a surface of signature $(0,+,4)$. So we get that the
action of $\sigma^*$ on $\tau_2$ is $\tau_2 \mapsto
-\overline{\tau_2}$. We have that the mapping $R$ lifts to
$A_2^{1/2}$, and the lift of $r$ is
$$\tilde{r}(z) = \frac{(1-\mu)\overline{z}+\mu}
{\mu\overline{z}+1+\mu}.$$ Computing as in the previous example, and
taking care of the fact that $r$ reverses orientation, we see that
$\sigma^*(\tau_1)$ is given by the cross ratio of the points
$$\begin{array}{l}
f(\tilde{r}A_2^{1/2} A_2^{-1}A_3^{-1} A_2^{-1/2}\overline{r}^{-1}),~
f(\tilde{r}A_2^{1/2} A_3 A_2^{-1/2}\overline{r}^{-1}),\\
f(A_2),~
\tilde{C_3}(f(\tilde{r}A_2^{1/2} A_2^{-1}A_3^{-1} A_2^{-1/2})),
\overline{r}^{-1})),\end{array}$$
where
$\tilde{C_3} = \tilde{r}A_2^{1/2} C_3^{-1} A_2^{-1/2}\overline{r}^{-1}$.
This cross ratio gives
$\sigma^*(\tau_1) = 1-\overline{\tau_3}$. Similarly, one gets
$\sigma^*(\tau_3) = 1-\overline{\tau_1}$. Therefore, the Teichm\"{u}ller space
$T(\Sigma)$ can be identified with the set of points
$(\tau_1,\tau_2,\tau_3) \in T(2,+,0)$ such that
$$\left\{\begin{array}{lcl}
{\rm Re}(\tau_2)&=&0\\
{\rm Re}(\tau_1)&=&1-{\rm Re}(\tau_3)\\
{\rm Im}(\tau_1)&=&{\rm Im}(\tau_3),\end{array}\right .$$
which proves theorem $1$ of the introduction.

The above computations give us some other isomorphisms, different from
those of the previous section. Observe that
the transformation $R$ is just the hyperelliptic involution on $\Sigma^c$.
It is not hard to see that $R^*$ acts like the identity in
$T(\Sigma^c)$
(\cite[pg. 126]{nag:teic}). The mapping $r$ has a curve of the partition
as the set of fixed points. We have that $\Sigma^c/<r>$ is
a sphere with one hole and two punctures, in the first example, or a
torus with a hole in the second example. Let us denote this surfaces by
$S_1$ and $S_2$, respectively.
\begin{cor}The spaces $T(S_1)$ and $T(S_2)$ are isomorphic to
$T(0,-,2)$ and $T(3,-,0)$ respectively.\end{cor}
Isomorphisms between deformation spaces of orientable and non-orientable
surfaces, as those of the above corollary, do not happen if the genus is
bigeer than $2$ (\cite[pg. 152]{sep:book}).

\ifx\undefined\bysame
\newcommand{\bysame}{\leavevmode\hbox to3em{\hrulefill}\thinspace}
\fi

Address: School of Maths, Tata Institute of Fundamental Research,
Bombay 4005, India

email: pablo@math.tifr.res.in
\nopagebreak\begin{figure}[ht]
\centerline{\psfig{figure=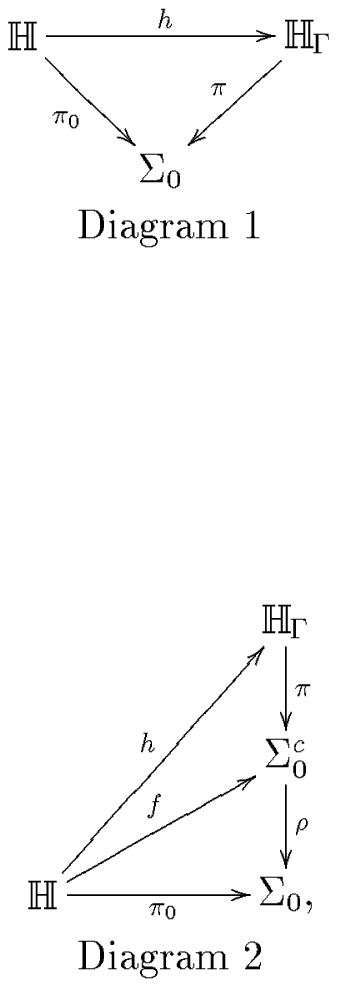,height=9in,width=6in}}
\end{figure}
\newpage
\begin{figure}[ht]
\centerline{\psfig{figure=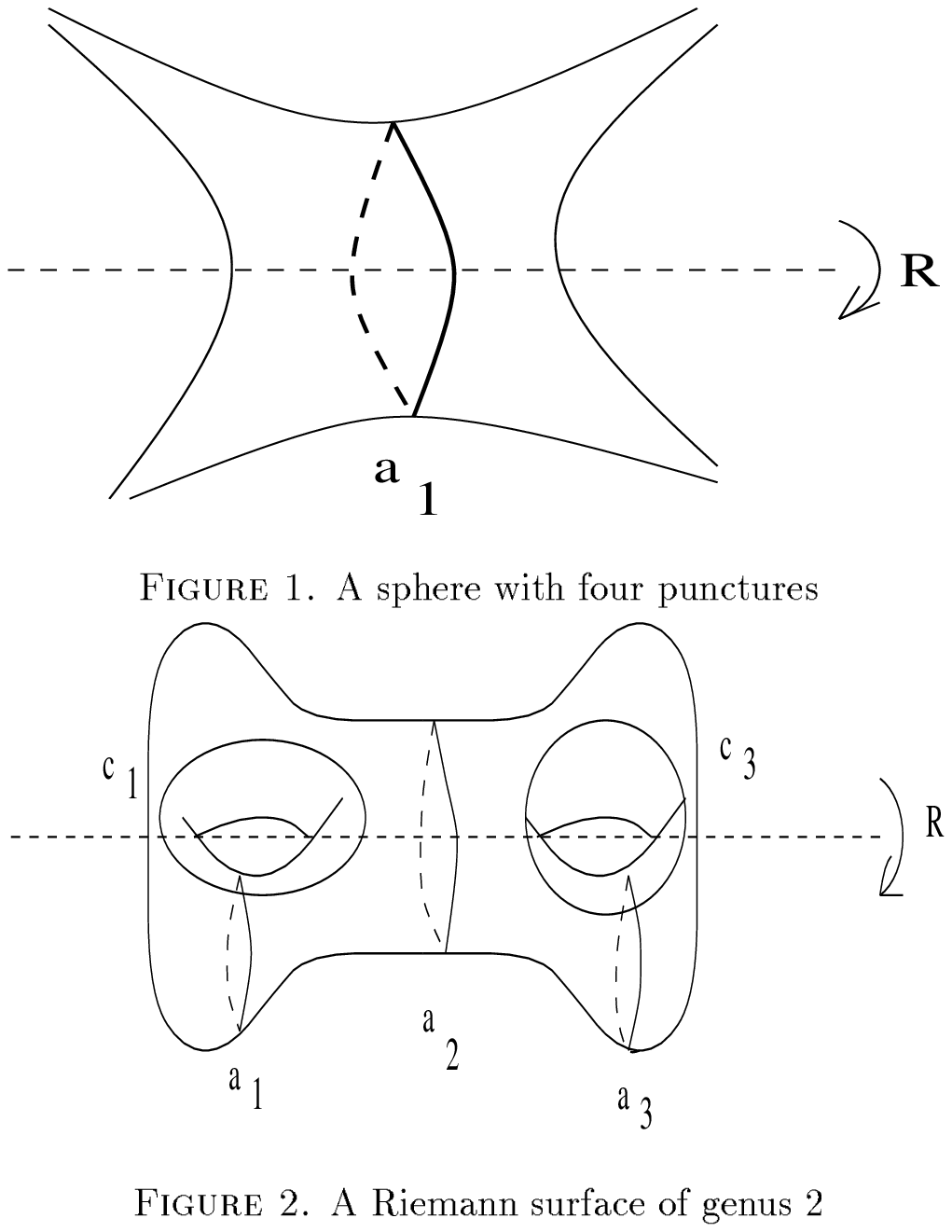,height=9in,width=6in}}
\end{figure}
\end{document}